\documentclass[11pt,reqno]{amsart}
\usepackage[utf8]{inputenc}
\usepackage[T1]{fontenc}
\usepackage{lmodern}
\usepackage[english]{babel}
\usepackage{amsmath,a4wide}
\usepackage{xfrac}
\usepackage{esint}
\usepackage{stmaryrd,mathrsfs,bm,amsthm,mathtools,yfonts,amssymb,color,braket,booktabs,graphicx,graphics,amsfonts}
\usepackage{latexsym,microtype,indentfirst,hyperref}
\usepackage{xcolor}
\usepackage{courier}
\usepackage{constants}
\usepackage{tikz}
\usepackage{multicol,caption}

\newtheorem{theorem}{Theorem}[section]

\theoremstyle{definition}
\newtheorem{definition}[theorem]{Definition}
\newtheorem{remark}[theorem]{Remark}

\newtheorem{assumption}[theorem]{Assumption}

\numberwithin{equation}{section}
\numberwithin{subsection}{section}

\newcommand{\R}{\mathbb{R}} 

\newcommand{\spt}{\mathrm{spt}} 
\newcommand{\dist}{\mathrm{dist}} 
\newcommand{\Ha}{\mathcal{H}} 

\newcommand{\mres}{\mathbin{\vrule height 1.6ex depth 0pt width 
0.13ex\vrule height 0.13ex depth 0pt width 1.3ex}}

\newcommand{\abs}[1]{\lvert#1\rvert} 

\newcommand{\var}{\mathbf{var}} 

\newconstantfamily{eps}{
symbol=\varepsilon}
\newconstantfamily{con}{
symbol=c}
\newconstantfamily{Lam}{
symbol=\Lambda}

\newenvironment{Figure}
  {\par\medskip\noindent\minipage{\linewidth}}
  {\endminipage\par\medskip}

\title[Existence and regularity of Brakke's mean curvature flows]{New advances on the existence and regularity of Brakke's mean curvature flows}
\date{\today}

\author[S. Stuvard]{Salvatore Stuvard}
\address{Dipartimento di Matematica, Universit\`{a} degli Studi di Milano, Via Saldini 50, I-20133 Milano (MI), Italy}
\email{salvatore.stuvard@unimi.it}

\begin{document}

\begin{abstract}
In this survey paper, I discuss some recent progress on the existence and regularity of Brakke flows. These include: an ``end-time version'' of Brakke's local regularity theorem, which allows to extend the validity of the celebrated regularity theorem by White from limits of smooth mean curvature flows to arbitrary Brakke flows; a global-in-time existence theorem for multi-phase Brakke flows of grain boundaries satisfying suitable ${\rm BV}$ regularity in time, in both the unconstrained and the fixed boundary settings, with applications to Plateau's problem for the latter; and the proof that branching singularities of minimal surfaces are a trigger for dynamical instability, in the sense that they may be ``perturbed away'' by a non-trivial canonical Brakke flow. This note is an extended version of a talk given by the author at the MATRIX Research Institute on the occasion of the workshop entitled ``\textit{Minimal surfaces and geometric flows: interaction between the local and the nonlocal worlds}''.
\end{abstract}

\maketitle

\section{Introduction}
Arising as the $L^2$-gradient flow of the area functional, the Mean Curvature Flow (henceforth often abbreviated MCF) is one of the most fundamental among the geometric flows involving extrinsic curvatures. Classically, a MCF of dimension $n \geq 1$ in a (smooth) Riemannian manifold $(\mathcal N, h)$ of dimension $n+k$ ($k \geq 1$) is a family of smoothly immersed $n$-dimensional submanifolds $M_t$, parametrized by time $t \in I$ with $I \subset \R$ an interval, such that the velocity of motion is equal to the mean curvature vector of the immersion at every point, for every time. It is a classical fact that, if $M_0$ is a smooth immersed submanifold of $\mathcal N$ then there exists a unique smooth solution to MCF starting with $M_0$ for short time. In fact, it is not difficult to construct examples of initial data such that the corresponding smooth MCF develops singularities in finite time, and before it becomes extinct: at such singular points, surfaces may undergo topology changes, so that any description in terms of the classical PDE ceases to make sense. The need to describe the motion by mean curvature of surfaces \emph{through} singularities and topology changes is the first motivation for studying \emph{weak solutions} to MCF.

\smallskip

\begin{multicols}{2}
{\begin{Figure}
\begin{center}
\vspace{1.1cm}
 \includegraphics[scale=0.45]{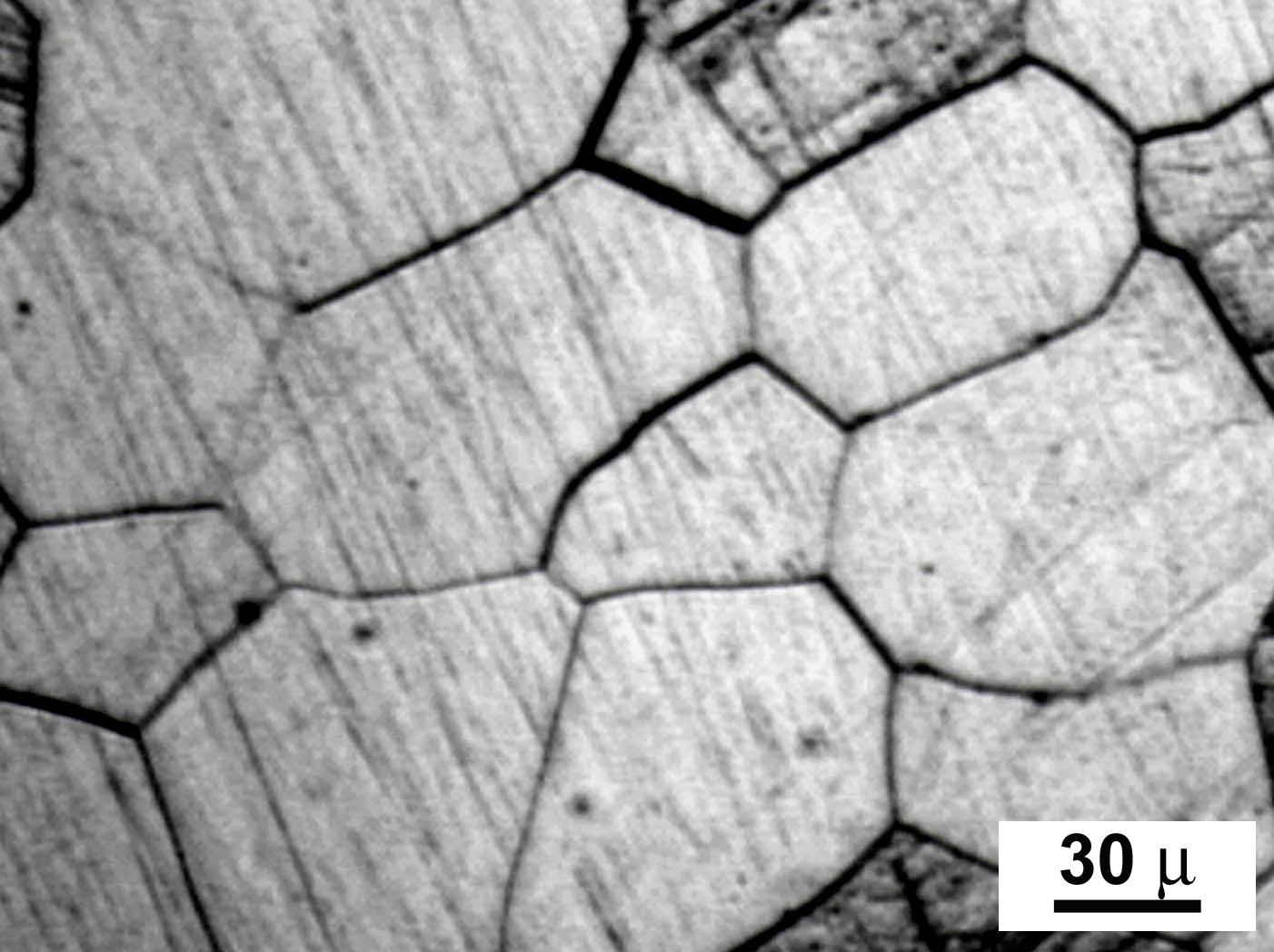}
 \captionof{figure}{\small{Grain boundaries in a metal. By Edward Pleshakov \href{https://creativecommons.org/licenses/by/3.0/deed.en}{CC BY 3.0}, Source \href{https://commons.wikimedia.org/w/index.php?curid=3912586}{Wikipedia}.}}\label{fig:grain-boundaries}
 \end{center}
\end{Figure}}
\columnbreak
Another motivation comes from the applications. The mean curvature flow has been studied for decades as a model for physical systems or materials governed by energies dominated by surface tension. Mullins may have been the first to write the MCF equation in \cite{Mullins}, when describing the motion of \emph{grain boundaries} and coarsening phenomena in metals; see Figure \ref{fig:grain-boundaries}. In these examples, the geometric objects of interest are intrinsically non-smooth. 
\end{multicols}

Brakke proposed a definition of weak solutions to MCF in his Ph.D. thesis \cite{Brakke_mcf}, under the scientific direction of Almgren. The definition is measure-theoretic in spirit: evolving immersed submanifolds are replaced by \emph{integral varifolds}, and in place of the PDE one finds a variational \emph{inequality}, which is in fact equivalent to the PDE if the varifolds are associated with smooth submanifolds. Such weak solutions are typically referred to as \emph{Brakke flows}.

\smallskip

This brief note is a review on some recent advances in the theory concerning existence and regularity of Brakke flows. 

\medskip

\noindent\textbf{Acknowledgements.} I would like to thank the MATRIX Research Institute for hospitality on the occasion of the workshop entitled ``\textit{Minimal surfaces and geometric flows: interaction between the local and the nonlocal worlds}''. My participation was supported through a MATRIX-Simons Travel Grant: I am thankful to the MATRIX Research Institute and the Simons Foundation. My work is supported by the \textit{Gruppo Nazionale per l'Analisi Matematica, la Probabilit\`a e le loro Applicazioni} of INdAM and by the project PRIN 2022PJ9EFL ``Geometric Measure Theory: Structure of Singular Measures, Regularity Theory and Applications in the Calculus of Variations''. 

\section{Preliminaries and plan of the paper} \label{s:def}

Let us assume that $\{M_t\}_{t \in [a,b)}$ is a one-parameter family of smooth immersed submanifolds of Euclidean space $\R^{n+k}$. Set $M_t = F_t (M)$, where $M$ is a smooth $n$-dimensional manifold, suppose that the immersion $F_t$ is the slice, at time $t \in [a,b)$, of a smooth map $F \colon M \times [a,b) \to \R^{n+k}$, and call $ v := \partial_t F$ the velocity of the flow. Let $U \subset \R^{n+k}$ be an arbitrary open set. The idea of Brakke is that the information that $v$ is the velocity of the flow is encoded in the equation which characterizes the rate of change of integrals of the form
\[
\int_{M_t} \phi(x,t) \, d\mathcal H^n(x)\,,
\]
where $\mathcal H^n$ denotes $n$-dimensional Hausdorff measure, and $\phi \in C^1(U \times [a,b))$ is such that $\phi(\cdot,t) \in C^1_{{\rm c}}(U)$. A calculation gives that
\begin{equation} \label{e:Brakke1}
    \frac{{d}}{dt} \int_{M_t} \phi(x,t) \, d\mathcal H^n(x) = \int_{M_t} \left\lbrace \left(\nabla\phi(x,t) - \phi(x,t)\,H(x,M_t)\right) \cdot v(x,t)  + \frac{\partial \phi}{\partial t}(x,t) \right\rbrace \,d\mathcal H^n(x) \,,
\end{equation}
where $H(x,M_t)$ is the mean curvature vector of $M_t$ at $x \in M_t$. In particular, for every $a \leq t_1 < t_2 < b$ it holds
\begin{equation}\label{e:Brakke2}
    \begin{split}
\int_{M_{t_2}} \phi(x,t_2) \, d\mathcal H^n  &- \int_{M_{t_1}} \phi(x,t_1) \, d\mathcal H^n \\
& \qquad  = \int_{t_1}^{t_2}\int_{M_t} \left\lbrace \left(\nabla \phi - \phi\,H(\cdot,M_t)\right) \cdot v  + \frac{\partial \phi}{\partial t} \right\rbrace \,d\mathcal H^n dt\,,
    \end{split}
\end{equation}
and if $M_t$ is a smooth MCF then \eqref{e:Brakke2} is satisfied for arbitrary $U$ and $\phi$ as above with $v$ replaced by $H(\cdot,M_t)$ in the integral on the right-hand side. 

\smallskip

Now, the interesting fact is that the identity in \eqref{e:Brakke2} can actually be translated so that it makes sense in the class of $n$-dimensional integral varifolds with generalized mean curvature vector $H \in L^2$. Let us recall that a \emph{varifold} of dimension $n$ in an open set $U \subset \R^{n+k}$ is a measure $V$ on the Grassmannian bundle $\mathbf{G}_n(U) := U \times {\rm G}(n+k, n)$, where ${\rm G}(n+k,n)$ is the Grassmannian of (unoriented) $n$-dimensional linear subspaces of $\R^{n+k}$. The varifold $V$ is \emph{integral} if the following holds: there are a countably $n$-dimensional rectifiable set $M \subset \R^{n+k}$ with locally finite $\mathcal H^n$-measure in $U$ and a \emph{multiplicity} function $\theta \in L^1_{{\rm loc}}(\mathcal H^n \mres_{M})$ with values in $\mathbb{Z}_{>0}$ such that
\begin{equation} \label{e:int_var}
    V(\varphi) = \int_M \varphi(x, {\rm Tan}(M,x)) \, \theta(x) \, d\mathcal H^n (x) \qquad \mbox{for every $\varphi \in C_{\rm c}({\bf G}_n(U))$}\,.
\end{equation}
In formula \eqref{e:int_var}, ${\rm Tan}(M,x)$ denotes the approximate tangent plane to $M$ at $x$, which exists for $\mathcal H^n$-a.e. $x \in M$ due to the fact that $M$ is rectifiable. We will write $V=\var (M,\theta)$ if $V$ is a varifold as in \eqref{e:int_var}. The \emph{first variation} of $V=\var(M,\theta)$ along a vector field $X \in C^1_{\rm c}(U ; \R^{n+k})$ is given by
\begin{equation}\label{e:fist_var}
    \delta V (X) = \int {\rm div}_{M} X \, d\|V\|\,.
\end{equation}
In \eqref{e:fist_var}, $\|V\|$ is the \emph{mass measure} $\|V\| = \theta\,\mathcal H^n \mres_{M}$ and ${\rm div}_M X$ is the tangential divergence of $X$ along $M$, defined for $\mathcal H^n$-a.e. $x \in M$ by
\[
{\rm div}_M X(x) : = {\rm trace}(DX^{\rm t} \circ {\rm Tan}(M,x))\,,
\]
where $DX$ is the Jacobian matrix of $X$, $DX^{{\rm t}}$ is its transpose, and ${\rm Tan}(M,x)$ denotes, by a slight abuse of notation, the projection operator $\R^{n+k} \to {\rm Tan}(M,x)$.

When the first variation $\delta V$ extends to a linear and continuous operator on $C_{{\rm c}}(U;\R^{n+k})$, one says that $V$ has \emph{locally bounded first variation}. By Riesz's representation theorem, $\delta V$ is then a $\R^{n+k}$-valued measure on $U$; if $\delta V$ is also absolutely continuous with respect to $\|V\|$, then we define the \emph{generalized mean curvature vector} of $V$ to be the Radon-Nikodym derivative
\begin{equation} \label{e:mean curvature}
    H(\cdot,V) = - \frac{d (\delta V)}{d\|V\|}\,.
\end{equation}

We can now give the definition of Brakke flow.

\begin{definition} \label{def:BF}
    Let $I = [a,b)$ be an interval, and let $U \subset \R^{n+k}$ be an open set. An \emph{$n$-dimensional integral Brakke flow} in $U$ is a one-parameter family of varifolds $\mathscr{V}=\{V_t\}_{t\in I}$ in $U$ such that all of the following hold: 
\begin{enumerate}
\item[(a)]
For a.e.~$t\in I$, $V_t$ is integral, $V_t = \var(M_t,\theta_t)$;
\item[(b)] For a.e.~$t\in I$, $\delta V_t$ is locally bounded and absolutely continuous with respect to $\|V_t\|$;
\item[(c)] The generalized mean curvature $H(\cdot,V_t)$ (which exists for a.e.~$t$ by (b)) satisfies $H(\cdot,V_t) \in L^2_{{\rm loc}}(\|V_t\|\times dt;\R^{n+k})$, and for every compact set $K \subset U$ and for every $t < b$ it holds $\sup_{s \in \left[a,t\right]} \|V_s\| (K) < \infty$;
\item[(d)]
For all $a\leq t_1<t_2<b$ and $\phi\in C_c^1(U\times\left[a,b\right);\mathbb R^+)$, it holds
\begin{equation}
\label{brakineq}
\begin{split}
&\|V_{t_2}\|(\phi(\cdot,t_2)) - \|V_{t_1}\|(\phi(\cdot,t_1)) \\ 
&\qquad\leq \int_{t_1}^{t_2} \int_{U} \left\lbrace \left( \nabla \phi(x,t) - \phi(x,t)\,H(x,V_t) \right) \cdot H(x,V_t) + \frac{\partial\phi}{\partial t}(x,t) \right\rbrace    d\|V_t\|(x)\,dt\,.
\end{split}
\end{equation}
\end{enumerate}
\end{definition}

The most striking difference between \eqref{e:Brakke2} and \eqref{brakineq} is that the latter is an \emph{inequality}, rather than an identity. In fact, it turns out that the inequality \eqref{brakineq} is \emph{equivalent} to the identity \eqref{e:Brakke2}, and thus to the classical PDE defining MCF, when $V_t$ are the unit multiplicity varifolds $V_t = \var(M_t,1)$ associated with \emph{smoothly immersed} $n$-dimensional submanifolds $M_t$. In the non-smooth setting, working only with the inequality is desirable for two reasons. One is purely technical: as the space-time integral of $\phi\,|H|^2$ (which appears with the minus sign on the right-hand side of \eqref{e:Brakke2} and \eqref{brakineq}) is only lower semicontinuous with respect to weak convergence, weak limits of smooth MCF are expected to only satisfy the inequality. The other motivation is more ``physical'', and it has to do with the fact that in a weak setting one may conceive systems evolving by mean curvature which exhibit singular regions disappearing at a much faster rate than the time-scale of the equation: at the level of the model, this disappearance will be regarded as instantaneous.

On one hand, it is remarkable that the inequality \eqref{brakineq} alone is sufficient to prove \emph{partial regularity results} for Brakke flows: a description of such results will be the content of Section \ref{s:regularity}. On the other hand, working only with the inequality has a fundamental drawback. Indeed, \eqref{brakineq} implies, among other things, that the difference in mass of $V_t$ between two instants of time $t_1 < t_2$ may be \emph{strictly less} than the amount dissipated in terms of the squared $L^2$-norm of the mean curvature over the interval $[t_1,t_2]$, and there is no control, a priori, on the extent to which instantaneous reduction of mass can occur. In particular, if $\mathscr{V}=\{V_t\}_{t \ge 0}$ is an integral Brakke flow then any flow $\tilde{\mathscr{V}}=\{\tilde V_t\}_{t \geq 0}$ defined by
\[
\tilde V_t = 
\begin{cases}
V_t & \mbox{for $0 \leq t < \tau$} \\
0 & \mbox{for $t \geq \tau$}\,,
\end{cases}
\]
where $\tau >0$ is arbitrary, is also an integral Brakke flow. In other words, Brakke's inequality alone permits \emph{sudden vanishing of the solution}, and is, therefore, a trigger for redundant non-uniqueness. Therefore, in practice it is desirable to work with classes of Brakke flows for which some sort of regularity in time can be guaranteed in order to avoid the occurrence of such phenomena. The approximation scheme proposed by Brakke in his general existence theorem of solutions to the Cauchy problem in \cite{Brakke_mcf} does not rule out the possibility that the presence of singularities in the initial datum triggers the catastrophic outcome described above at some point of the evolution. In selected geometric settings, particularly in the codimension $k=1$ case, some progress has been made by proposing either alternative notions of weak solution or alternative approximation schemes. A notion of \emph{viscosity solution} to MCF was proposed independently in \cite{ES} and \cite{CGG}: the methods therein produce a time-parametrized family of closed sets which are the \emph{level sets} (corresponding to a given level) of the \emph{unique} viscosity solution of a parabolic PDE, which is defined also after the occurrence of singularities. It is possible that the closed set may develop nontrivial interior afterwards, a phenomenon called \emph{fattening}, and it is not clear if the sets are a Brakke flow after singularities appear in general. On the other hand, Evans and Spruck proved in \cite{ES_IV} that, given the viscosity solution, for almost all levels the corresponding level sets are a unit multiplicity Brakke flow. On the side of devising alternative methods for constructing Brakke flows, it is imperative to mention \emph{elliptic regularization} \cite{Ilm1}, as well as \emph{phase field approximations} via the (parabolic) Allen–Cahn equation \cite{Ilm_AC,Ton_integrality}. While these methods produce global-in-time existence results of Brakke flows, they all rely on the ansatz that the MCF is represented as the boundary of a \emph{single} time-parametrized set, and therefore they do not allow the handle the framework of \emph{multi-phase} flows, which, on the other hand, is of fundamental importance for the applications, for instance in Materials Science, in the context of describing the dynamics of grain boundaries as in Figure \ref{fig:grain-boundaries}. Section \ref{s:canonical} will discuss some recent advances in the existence theory for multi-phase Brakke flows of grain boundaries, initiated by the work of Kim and Tonegawa; see \cite{kim-tone}. Finally, in Section \ref{s:dynamical} we will present results concerning the global-in-time existence of $n$-dimensional integral Brakke flows in an open and strictly convex domain $U \subset \R^{n+1}$ under the additional constraint that the evolving varifolds have \emph{fixed boundary} (defined in a suitable topological sense) on $\partial U$. The asymptotic analysis of the flow as $t \to \infty$ gives rise to a novel approach to the existence of varifold solutions to Plateau's problem for a large class of boundaries. This approach is purely dynamical, and it does not rely on classical minimization or min-max procedures. Furthermore, it naturally leads to introducing a notion of \emph{dynamical stability} for minimal surfaces, which appears to be powerful enough to prevent the formation of certain singularity types, and deserves further investigation.

\section{Partial regularity of unit multiplicity Brakke-like flows} \label{s:regularity}

A basic technique, common to a plethora of problems in Geometric Analysis, to investigate whether a solution is locally regular (that is, ``smooth'' or ``as smooth as classical solutions'') at a given point is commonly referred to as \emph{blow-up}: it consists in deducing information on the solution from the analysis of possible limits of sequences of suitable rescalings of the latter at the given point. In order to conclude the existence of blow-up limits, a suitable \emph{monotonicity formula} is typically needed: in the context of MCF, such monotonicity formula was proved by Huisken in \cite{Huisken_mono} and later extended to Brakke flows by Ilmanen in \cite{Ilmanen_mono}. For future reference, we record here its statement, working under the simplifying assumption that the Brakke flow $\mathscr{V}=\{V_t\}_{t \in [a,b)}$ in $\R^{n+k}$ satisfies the condition that, for all $t \in [a,b)$, $\spt\|V_t\|$ is contained in some bounded domain of $\R^{n+k}$, say $B_R$ for $R>0$ sufficiently large. For $x,y \in \R^{n+k}$ and $t < s$, we define the backwards $n$-dimensional heat kernel with pole $(y,s)$ by
\begin{equation} \label{e:heat_kernel}
    \varrho_{(y,s)}(x,t) := \frac{1}{(4\pi (s-t))^{\sfrac{n}{2}}} \,\exp\left(- \frac{|x-y|^2}{4 (s-t)} \right)\,.
\end{equation}

\begin{theorem}[Huisken's monotonicity formula] \label{thm:Huisken_mf}
    Let $\mathscr V = \{V_t\}_{t \in [a,b)}$ be an $n$-dimensional integral Brakke flow in $\R^{n+k}$ such that $\spt\|V_t\| \subset B_R$ for all $t \in [a,b)$. Let $y \in \R^{n+k}$ and $s \in (a,b]$. Then, for every $a \leq t_1 < t_2 < s$ it holds
    \begin{equation}\label{e:Huisken_mf}
        \begin{split}
            & \int \varrho_{(y,s)}(x,t_2) \, d\|V_{t_2}\|(x) - \int \varrho_{(y,s)}(x,t_1) \, d\|V_{t_1}\|(x) \\
            &\qquad \qquad \leq - \int_{t_1}^{t_2} \int \Big| H (x,V_t) + \frac{S^\perp(x-y)}{2(s-t)} \Big|^2 \, \varrho_{(y,s)}(x,t) \, dV_t(x,S)\,dt\,.
        \end{split}
    \end{equation}
In particular, the function $t \in [a,s) \mapsto \int \varrho_{(y,s)}(x,t)\, d\|V_t\|(x)$ is non-increasing. The limit
\begin{equation} \label{e:GD}
    \Theta(\mathscr{V}, (y,s)) := \lim_{t \to s^-} \int \varrho_{(y,s)}(x,t)\, d\|V_t\|(x)
\end{equation}
is called the \emph{Gaussian density} of $\mathscr{V}$ at $(y,s)$.
\end{theorem}

Suppose that $(y,s)$ belongs to the support of the product measure $\|V_t\|\times dt$ in $\R^{n+k} \times \R$, assume without loss of generality that $[a,b)=[-1,1)$ and $(y,s) = (0,0)$, and consider, for every $\lambda > 0$, the parabolic rescalings of $\mathscr V$ around $(0,0)$ with factor $\lambda$, that is the flows $\{V^{\lambda}_\tau\}_\tau$ defined by
\[
V^{\lambda}_\tau := (\eta_{\lambda})_\sharp V_{\lambda^2 \tau} \qquad \mbox{for $\tau \in [-\lambda^{-2},\lambda^{-2})$}\,,
\]
where $\eta_{\lambda}(x) := \lambda^{-1}x$ and $(\eta_{\lambda})_\sharp$ denotes the push-forward map on varifolds corresponding to $\eta_{\lambda}$, see \cite{Simon}. Thanks to Theorem \ref{thm:Huisken_mf}, for every sequence $\lambda_j \to 0$ the masses of $V^{\lambda_j}_\tau$ are equi-bounded uniformly in $\tau \in [-L,0)$ for every $L > 0$, and thus the compactness theorem of Brakke flows (see e.g. \cite[Theorem 3.7]{Ton1}) guarantees that a (not relabeled) subsequence of $\mathscr{V}^{\lambda_j} = \{V^{\lambda_j}_\tau\}_\tau$ converges, as $j \to \infty$, to a limit Brakke flow $\hat{\mathscr{V}} = \{\hat V_\tau\}_\tau$ defined and selfsimilar for $\tau \in (-\infty,0)$: $\hat{\mathscr{V}}$ is a \emph{tangent flow to $\mathscr V$ at $(0,0)$}. Notice that, \emph{a priori}, tangent flows may be non-unique, namely they may depend on the particular sequence $\lambda_j$. By the dimension reduction argument of White \cite{White_stratification}, at \emph{most} points $(y,s) \in \spt(\|V_t\|\times dt)$ the Brakke flow $\mathscr{V}$ admits a tangent flow which is static and an integer multiple of a plane: that is, $\hat V_\tau$ is defined for every $\tau \in \R$ and is identically equal to $\var(\pi,Q)$, where $Q \in \mathbb{Z}_{>0}$ and $\pi$ is an $n$-dimensional linear subspace of $\R^{n+k}$. The word ``most'' here can be intended in terms of \emph{parabolic Hausdorff dimension}, namely Hausdorff dimension in the metric space $(\R^{n+k} \times \R, {\rm d})$, where 
\begin{equation}\label{def:parabolic metric}
    {\rm d}((x_1,t_1), (x_2,t_2)) := \max\{|x_1-x_2| , \sqrt{|t_1-t_2|}\}\,.
\end{equation}
In fact, the set of points $(y,s) \in \spt(\|V_t\|\times dt)$ which do not admit static tangent planes has parabolic Hausdorff dimension at most $n+1$ (informally, the flow itself has parabolic dimension $n+2$,
as time is $2$-dimensional with respect to ${\rm d}$). Let us call \emph{regular} a point $(y,s) \in \spt(\|V_t\| \times dt)$ which has a neighborhood in space-time where $\mathscr{V}$ is a smooth MCF (or a constant integer multiple thereof). If $(y,s)$ is a regular point, then necessarily $\mathscr{V}$ has, at $(y,s)$, a unique tangent flow which is a static plane, possibly with constant multiplicity. It is then natural to ask if the converse is true: namely, if a point for which a blow-up is a static plane is necessarily regular. It turns out that this implication holds true if, and only if, the multiplicity of the blow-up is $Q=1$. Tangent planes with multiplicity $Q\ge 2$ may arise at singular points even when the Brakke flow is time-independent, namely $V_t \equiv V_0$ for a stationary integral varifold $V_0$. In codimension $k=2$, classical examples of this phenomenon (commonly referred to as \emph{branching}) are given by (the unit multiplicity varifolds associated with) holomorphic varieties such as $M_0 = \{(z,w) \in \mathbb C^2 \, \colon \, w^2=z^3\}$, which are even locally area minimizing in the sense of integral currents (once equipped with a suitable orientation); in codimension $k=1$, branching cannot occur in the setting of area minimizing integral currents, but there are examples in the context of currents minimizing the area in their $\mathbb{Z}_p$-homology class when $p$ is an even integer, which induce stable stationary varifolds; see, for instance, \cite[Example 1.6]{DLHMS}. If the multiplicity of the tangent plane is $Q=1$, regularity for stationary varifolds is a corollary of the celebrated $\varepsilon$-regularity theorem by Allard \cite{Allard}. In the context of Brakke flows, there is an analogue of Allard's theorem, known as \emph{Brakke's local regularity theorem}, Theorem \ref{thm:Brakke_reg} below. The theorem holds, in fact, for a larger class of ``Brakke-like'' flows, which solve a PDE of the form
\[
\mbox{velocity $=$ mean curvature $+$ forcing term}\,.
\]
More precisely, one can consider flows $\mathscr V = \{V_t\}_{t \in [a,b)}$ as in Definition \ref{def:BF} with condition ${\rm (d)}$ replaced by the following:
\begin{enumerate}
\item[(d)']  For all $a \leq t_1 < t_2 < b$ and $\phi \in C^1_{{\rm c}}(U \times [a,b) ; \R^+)$, it holds
\begin{equation}
\label{brakineq_forcing}
\begin{split}
&\|V_{t_2}\|(\phi(\cdot,t_2)) - \|V_{t_1}\|(\phi(\cdot,t_1)) \\ 
&\quad\leq \int_{t_1}^{t_2} \int_{U} \left\lbrace \left( \nabla \phi(x,t) - \phi(x,t)\,H(x,V_t) \right) \cdot \left(H(x,V_t) + S^\perp u(x,t) \right) + \frac{\partial\phi}{\partial t}(x,t) \right\rbrace    d V_t (x,S)\,dt\,,
\end{split}
\end{equation}
where $u \colon U \times [a,b) \to \R^{n+k}$ is a function in a Lebesgue class $L^{p,q}_{x,t}$, namely such that
\[
\|u\|_{L^{p,q}}:=\left(\int_a^b \left( \int_U |u(x,t)|^p \, d\|V_t\|(x) \right)^{\sfrac{q}{p}} \, dt\right)^{\sfrac{1}{q}} < \infty
\]
for $2 \leq p < \infty$ and $2 < q < \infty$ satisfying
\begin{equation}\label{pq_condition} 
\zeta := 1 - \frac{n}{p} - \frac{2}{q} > 0\,,
\end{equation}

or with the analogous condition on the natural norm in case $p$ or $q$ are infinite.
\end{enumerate}

The need to treat Brakke flows with forcing term arises naturally when one wants to consider MCF and the corresponding notion of Brakke flow in a general
$(n+k)$-dimensional Riemannian manifold $\mathcal N$. By Nash's 
isometric embedding theorem, we may always consider $\mathcal N$ to be
a submanifold in a domain $U\subset \mathbb R^d$ for some sufficiently large $d$. A Brakke flow in $\mathcal N$
can be then defined by asking ${\rm spt}\|V_t\|
\subset \mathcal N$ for all $t$, and by requiring the validity of \eqref{brakineq_forcing} with $u(x,t) = - \sum_{i=1}^{n} A(\tau_i,\tau_i)$, where $\{\tau_1,\ldots,\tau_n\}$ is an orthonormal basis of ${\rm Tan}(M_t,x)$ and $A$ is the second fundamental form of the immersion $\mathcal N \hookrightarrow \R^d$; see also \cite[Section 7]{Ton-2}.

\begin{theorem}[Brakke's local regularity theorem] \label{thm:Brakke_reg}
    For every $\nu \in (0,1)$ there exists $\varepsilon(n,k,p,q,\nu) \in (0,1)$ with the following property. Suppose $\mathscr{V}=\{V_t\}_{t \in [-1,0)}$ is a flow as in Definition \ref{def:BF} with $I=[-1,0)$, $U = B_2^n(0)\times \R^k \subset \R^n \times \R^k \simeq \R^{n+k}$, and with ${\rm (d)}$ replaced by ${\rm (d)'}$. Assume also that $\mathscr V$ is \emph{unit multiplicity}, namely that $V_t = \var(M_t,1)$ for a.e. $t \in I$. Denoting $\pi_0 := \R^n \times \{0\}$ and ${\rm C}_r(0,\pi_0)$ the cylinder $B_r^n(0)\times \R^k$, suppose the following:
\begin{itemize}
\item[(H1)] it holds
\[
\mu^2 := \int_{-1}^0 \int_{{\rm C}_2(0,\pi_0)} \dist^2(x,\pi_0)\,d\|V_t\|(x)\,dt < \varepsilon^2\,, \qquad \|u\|_{L^{p,q}} < \varepsilon\,; 
\]
\item[(H2)] $\|V_{-4/5}\|({\rm C}_1(0,\pi_0)) \leq (2-\nu) \omega_n$, where $\omega_n = \mathcal L^n(B_1)$;
\item[(H3)] $\spt(\|V_t\| \times dt) \cap ({\rm C}_{\nu}(0,\pi_0)\times\{0\}) \neq \emptyset$.
\end{itemize}
Then, for all $t \in [-\sfrac14,0)$, $\spt\|V_t\| \cap {\rm C}_{\sfrac12}(0,\pi_0)$ coincides with the graph of a $C^{1,\zeta}$ function $f=f(z,t)$ defined on $B^n_{\sfrac12}(0) \times [-\sfrac14, 0)$ with $\|f\|_{C^{1,\zeta}} \lesssim \mu$.
\end{theorem}

\begin{remark}
    If $u$ belongs to some H\"older class $C^{\ell,\alpha}$, then one obtains $C^{\ell+2,\alpha}$ estimates on $f$ assuming smallness of the norm $\|u\|_{C^{\ell,\alpha}}$ in ${\rm (H1)}$. If $u \equiv 0$ then $f$ is smooth with all derivatives bounded in terms of $\mu$, and it solves MCF classically.
\end{remark}

A statement of Theorem \ref{thm:Brakke_reg} is contained in Brakke's work \cite{Brakke_mcf} for the case $u \equiv 0$. Since Huisken's monotonicity formula was not available at the time, Brakke's proof relies on a chain of clever graphical approximations of the flow, and it may be hard to follow. A revisited proof following the classical De Giorgi-type blow-up method also employed by Allard in \cite{Allard} in the setting of stationary varifolds, and with a non-zero forcing $u$, was given by Kasai and Tonegawa in \cite{Kasai-Tone}. Nonetheless, under the same assumptions of Theorem \ref{thm:Brakke_reg}, Kasai and Tonegawa obtained graphicality only in a smaller time interval, say for $t \in [-\sfrac{1}{4}, - \sfrac{1}{8}]$, and thus with some ``waiting time'' for regularity towards the end of the interval where the flow is defined. The regularity up to the end-time was proved by Tonegawa and the author in \cite{ST_endtime}. An alternative proof up to $C^{1,\zeta}$ regularity in the case of forcing term $u \in L^\infty$ was given by De Philippis, Gasparetto, and Schulze in \cite{DPGS}, using the viscosity techniques introduced by Savin in the context of elliptic PDE. Notice that the hypothesis ${\rm (H3)}$ is needed to guarantee that the flow does not vanish up until the end-time.

The proof of the end-time regularity allowed Tonegawa and the author to prove, in \cite{ST_endtime}, the following theorem, which extends the celebrated local regularity theorem by White in \cite{White_regMCF} to arbitrary Brakke flows, and not only to those which arise as limits of smooth mean curvature flows.

\begin{theorem}[White's local regularity theorem] \label{thm:White_reg}
   There exists $\varepsilon=\varepsilon(n,k) > 0$ with the following property. If $\mathscr{V}=\{V_t\}_t$ is a unit multiplicity integral Brakke flow of dimension $n$ in an open domain $U$ of Euclidean space $\R^{n+k}$ or a Riemannian manifold of dimension $(n+k)$, if $(y,s)\in \spt(\|V_t\|\times dt)$, and if $\Theta(\mathscr{V},(y,s)) \leq 1 + \varepsilon$, then there exists $r > 0$ such that $\mathscr{V}$ is a smooth mean curvatur flow in $B_r(y)\times (s-r^2,s)$.
\end{theorem}

The following partial regularity theorem of unit multiplicity Brakke flows is a consequence of Brakke's local regularity theorem; see \cite{Kasai-Tone}.

\begin{theorem}[Partial regularity] \label{thm:partial_reg}
    If $\mathscr{V}=\{V_t\}_{t \in I}$ is an $n$-dimensional Brakke-like flow in an open subset $U \subset \R^{n+k}$ with forcing $u \in L^{p,q}_{x,t}$ (with $p,q$ satisfying \eqref{pq_condition}) then for a.e. $t \in I$ there exists a (possibly empty) closed set $G_t\subset \spt\|V_t\|$ with $\mathcal H^n(G_t)=0$ such the point $(x,t)$ is regular for $\mathscr{V}$ for every $x \in \spt\|V_t\|\setminus G_t$.
\end{theorem}

\section{Existence of canonical multi-phase Brakke flows} \label{s:canonical}

In this section, we present some recent results on the global-in-time existence of multi-phase flows of grain boundaries. Compared to the previous sections, we immediately set the codimension to be $k=1$. The initial data we can consider are essentially arbitrary $n$-rectifiable sets in $\R^{n+1}$. More precisely, we will work under the following assumption.

\begin{assumption}\label{ass:existence}
    We assume that:
    \begin{itemize}
\item[(A1)] $M_0 \subset \R^{n+1}$ is a relatively closed, countably $n$-rectifiable set, and $\Ha^n (M_0 \cap B_R)$ grows at most exponentially fast as $R \to \infty$: more precisely, there exists a function $\Omega \in C^2(\R^{n+1})$ and a constant $c \geq 0$ with
\[
0 < \Omega(x) \leq 1\,, \quad |\nabla \Omega(x)| \leq c\,\Omega(x)\,, \quad \|\nabla^2 \Omega(x)\|\leq c\,\Omega(x)\qquad \mbox{for every $x\in \R^{n+1}$}
\]
such that
\[
\mathcal H^n \mres_{\Omega}(M_0) := \int_{M_0} \Omega(x)\,d\mathcal H^n(x) < \infty\,;
\]
\item[(A2)] $\R^{n+1} \setminus M_0 = \bigcup_{i=1}^N E_{i,0}$, where $E_{i,0}$ are non-empty, open and mutually disjoint subsets of $\R^{n+1}$, and $N \geq 2$.
\end{itemize}
\end{assumption}

The following existence theorem is contained in the works of Kim and Tonegawa \cite{kim-tone} and Tonegawa and the author \cite{ST_canonical}.

\begin{theorem}[Existence of canonical multi-phase Brakke flows] \label{thm:canonical}
    Under Assumption \ref{ass:existence}, there exist an $n$-dimensional integral Brakke flow $\mathscr{V}=\{V_t\}_{t \ge 0}$ in $\R^{n+1}$ as well as flows $\{E_i(t)\}_{t \geq 0}$ for $i \in \{1,\ldots,N\}$ of open sets with the following properties. Let $\sigma$ denote the product measure $\sigma = \|V_t\| \times dt$ in $\R^{n+1}\times \R$, and $M_t := \R^{n+1} \setminus \bigcup_{i=1}^N E_i(t)$. Then:
    \begin{itemize}
        \item[(i)] \textbf{initial condition:} $V_0 = \var(M_0,1)$, and $E_i(0) = E_{i,0}$ for every $i$;
        \item[(ii)] \textbf{continuity of mass at the initial time:} if $\mathcal{H}^n(\bigcup_i (\partial E_{i,0} \setminus \partial^*E_{i,0})) = 0$, where $\partial^*$ denotes reduced boundary, then $\lim_{t \to 0^+} \|V_t\| = \|V_0\| = \mathcal H^n \mres_{M_0}$;
        \item[(iii)] \textbf{mass and total mean curvature bounds}: $\|V_t\|(\Omega) \leq \mathcal H^n \mres_{\Omega}(M_0) \, \exp(\sfrac{c^2 t}{2})$, and
        \[
        \int_0^t\int_{\R^{n+1}} |H(x,V_s)|^2\,\Omega(x)\,d\|V_s\|(x)\,ds < \infty \qquad \mbox{for every $t > 0$}\,;
        \]
        \item[(iv)] \textbf{dissipation inequality:} if $\mathcal H^n(M_0)<\infty$, and thus the choice $\Omega \equiv 1$ is possible in Assumption \ref{ass:existence}, then
        \[
        \|V_t\|(\R^{n+1}) \leq \mathcal H^n(M_0) -  \int_0^t\int_{\R^{n+1}} |H(x,V_s)|^2\,d\|V_s\|(x)\,ds \qquad \mbox{for every $t > 0$}\,;
        \]
        \item[(v)] \textbf{grain boundaries:} for every $t > 0$, $E_1(t),\ldots, E_N(t)$ are pairwise disjoint; moreover, \[M_t = \bigcup_{i} \partial E_i(t) = \{x \in \R^{n+1} \, \colon \, (x,t) \in \spt(\sigma)\}\,;\]
        \item[(vi)] \textbf{support of the flow:} for every $t > 0$, $\spt\|V_t\| \subset M_t$ and $\mathcal H^n(M_t \cap K) < \infty$ for every $K \subset \R^{n+1}$ compact; moreover, for a.e. $t > 0$ $\spt\|V_t\|$ is countably $n$-rectifiable and $\mathcal H^{n-1+\delta}(M_t \setminus \spt\|V_t\|) = 0$ for every $\delta > 0$, so that \[V_t=\var(M_t,\theta_t) \quad \mbox{with $\theta_t(x) = \Theta^n(\|V_t\|,x)$ at $\|V_t\|$-a.e. $x$} \qquad \mbox{for a.e. $t > 0$}\,\]
        where $\Theta^n(\|V_t\|,x) := \lim_{r \to 0^+} (\omega_n r^n)^{-1} \,\|V_t\|(B_r(x))$ is the $n$-dimenional density of the measure $\|V_t\|$ at $x$;
        \item[(vii)] \textbf{generalized ${\rm BV}$-flow property:} for every $i \in \{1,\ldots,N\}$, for every $0 \leq t_1 < t_2 < \infty$, and for every $\phi \in C^1_{\rm c}(\R^{n+1} \times [0,\infty))$ it holds
        \begin{equation} \label{e:canonical}
            \begin{split}
               & \int_{E_i(t_2)} \phi(x,t_2) \, dx - \int_{E_i(t_1)} \phi(x,t_1) \, dx \\
                & \qquad = \int_{t_1}^{t_2} \left( \int_{E_i(t)} \frac{\partial\phi}{\partial t} (x,t) \, dx + \int_{\partial^*E_i(t)} \phi(x,t)\, H(x,V_t) \cdot \nu_{E_i(t)}(x) \, d\mathcal{H}^n(x) \right) \, dt\,,
            \end{split}
        \end{equation}
        where $\nu_{E_i(t)}(x)$ is the outer unit normal vector to $\partial^*E_i(t)$ at $x$;
        \item[(viii)] \textbf{two-sidedness at unit density points:} if $N \geq 3$, then for a.e. $t > 0$ and $\mathcal H^n$-a.e. $x\in M_t$ the condition $\theta_t(x)=1$ implies $x \in \bigcup_i \partial^*E_i(t)$;
        \item[(ix)] \textbf{two-phase case:} if $N=2$, then for a.e. $t > 0$ and $\mathcal H^n$-a.e. $x \in M_t$ it holds
        \[
        \theta_t(x) = \begin{cases}
            \mbox{odd integer} & \mbox{for $x \in \partial^*E_1(t)\, (=\partial^*E_2(t))$}, \\
            \mbox{even integer} & \mbox{for $x \in M_t \setminus \partial^*E_1(t)$}\,.
        \end{cases}
        \]
    \end{itemize}
\end{theorem}

We remark that the identity \eqref{e:canonical} implies, in particular, that for every ball $B\subset\R^{n+1}$ and for every $i \in \{1,\ldots,N\}$ the map $t \in (0,\infty) \mapsto \mathcal L^{n+1}(E_i(t)\cap B)$ is $\sfrac12$-H\"older continuous, so that the Brakke flows of Theorem \ref{thm:canonical} cannot experience sudden vanishing of all the mass. Rather, grains evolve continuously in time, and their reduced boundary is advected precisely by the generalized mean curvature vector of the underlying varifolds. In fact, for every $i$ the function $\chi_i(x,t) = \chi_{E_i(t)}(x)$ (where $\chi_E$ denotes characteristic function of $E$) is ${\rm BV}_{{\rm loc}}(\R^{n+1}\times [0,\infty))$, with space-time gradient
\[
D\chi_i(x,t) = (-\nu_{E_i(t)}(x), H(x,V_t) \cdot \nu_{E_i(t)}(x))\, d\mathcal H^n \mres_{\partial^*E_i(t)}\,dt\,.
\]
In the two-phase case ($N=2$), existence of a generalized velocity $v \in L^2$ moving the common boundary was proved by Mugnai-R\"oger for Brakke flows arising as sharp interface limits of minimizers of the Allen-Cahn action functional in \cite{MugRog}; more recently, Hensel-Laux have shown that the ${\rm BV}$ property with velocity in $L^2$ holds for limits of solutions to the Allen-Cahn equation in \cite{HenselLaux_ACisBV}. In any case, even in the two-phase case, the characterization $\partial_t \chi_i= (H \cdot \nu_{E_i}) d\mathcal H^n\mres_{\partial^*E_i}dt$ was not known in general.

\smallskip

The Brakke flow $\mathscr V = \{V_t\}_{t \geq 0}$ and the flow of grains $\{E_i(t)\}_{t \geq 0}$ of Theorem \ref{thm:canonical} are obtained by means of an algorithm, devised by Kim and Tonegawa in \cite{kim-tone}, which produces, starting with the initial datum $(M_0, \{E_{i,0}\}_{i=1}^N)$, a sequence of approximate solutions to MCF. More precisely, for each sufficiently large $j$ (the index of the sequence), the algorithm produces a flow $\{\mathcal E_j (t)\}_{t \in \left[0,j\right]}$, where $\mathcal E_j(t)$ is an $\mathcal L^{n+1}$-partition of $\R^{n+1}$ in open sets, namely $\mathcal E_j(t) = \{E_{j,1}(t),\ldots, E_{j,N}(t)\}$. Moreover, partitions are piecewise constant in time, namely $\mathcal E_{j}(t) = \mathcal E_{j}^k = \{E_{j,i}^k\}_{i=1}^N$ for $t$ in intervals (\emph{epochs}) $\left( (k-1)\, \Delta t_j, k\,\Delta t_j \right]$ of length $\Delta t_j \to 0^+$. The idea is then that if the partition $\mathcal E_{j}^{k+1}$ is constructed (inductively) from the partition $\mathcal E_{j}^k$ appropriately then, along a suitable subsequence $j_\ell \to \infty$, the (varifolds associated to) the boundaries $\partial \mathcal E_{j_\ell}(t)$ converge to the desired Brakke flow $V_t$. The open partition $\mathcal E_{j}^{k+1}$ at a given epoch is constructed from the open partition $\mathcal E_{j}^k$ at the previous epoch by applying two operations, which we call \emph{steps}. The first step is a small \emph{Lipschitz deformation} of partitions with the effect of regularizing singularities by locally minimizing the area of the boundary of partitions at a suitably small scale; the second step consists of flowing the boundary of partitions by an appropriately defined \emph{approximate mean curvature vector}, obtained by smoothing the surfaces' first variation via convolution with a localized heat kernel. The only difference between the scheme used in \cite{kim-tone} and the one adopted in \cite{ST_canonical} is in the class of admissible Lipschitz deformations in the first step: while in \cite{kim-tone} one only requires that the change of volumes of the grains due to Lipschitz deformation is small (for a certain smallness scale), in \cite{ST_canonical} we ask that the change of volume is small \emph{compared to the reduction in surface measure}. This is crucial in order to estimate that the contribution of the Lipschitz deformation step to the variation of integrals of test functions on the bulk of each grain along the approximation vanishes in the limit as $j_\ell \to \infty$ and prove \eqref{e:canonical}.

\smallskip

We conclude this section by remarking that it is tempting to conjecture that the regularization step by Lipschitz deformations embedded in the approximation scheme may have effects on the regularity properties of the Brakke flows produced in \cite{kim-tone} and \cite{ST_canonical}. A positive result in this direction is the following theorem of Kim and Tonegawa \cite{KiTo20}, valid for flows in dimension $n=1$.

\begin{theorem}[Improved regularity of network flows] \label{thm:network_regularity}
    Suppose $n=1$, and let $\mathscr{V} = \{V_t\}_{t \geq 0}$ be the Brakke flow obtained in \cite{kim-tone,ST_canonical} as in Theorem \ref{thm:canonical}. Then, for a.e. $t > 0$ the varifold $V_t$ has the following property. For every point $x \in \spt\|V_t\|$, there is a neighborhood $B_r(x)$ such that $V_t\mres_{B_r(x)}$ is the union of finitely many curves of class $W^{2,2}$ (and thus $C^{1,\sfrac12}$) which meet at $x$ forming angles of either $0$, $60$, or $120$ degrees.
\end{theorem}

\section{Brakke flows with fixed boundary, applications to Plateau's problem, and the notion of dynamical stability of minimal surfaces}\label{s:dynamical}

A natural question stemming from the theory detailed in Section \ref{s:canonical} is whether suitable modifications of the scheme in \cite{kim-tone,ST_canonical} may be used to produce existence results of Brakke-like flows under various types of constraints. In \cite{ST19}, Tonegawa and the author showed that this is indeed possible in the context of mean curvature flow whose boundary is kept fixed on the boundary of a strictly convex domain. 

\begin{figure}
 \centering
 \includegraphics[scale=0.35]{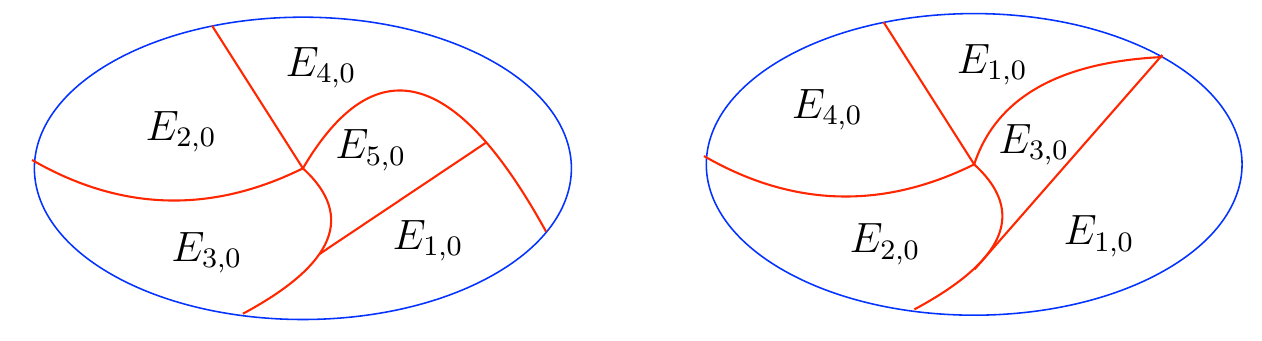}
 \caption{Assumption ${\rm (A3)}$: an admissible initial configuration on the left, and a not admissible one on the right.} \label{fig:bdry_cond}
\end{figure}

\begin{theorem}[Existence of Brakke flows with fixed boundary] \label{thm:Brakke_Dirichlet}
    Let $U \subset \R^{n+1}$ be a strictly convex, bounded open set with boundary $\partial U$ of class $C^2$. Assume that:
        \begin{itemize}
            \item[(A1)'] $M_0 \subset U$ is a relatively closed, countably $n$-rectifiable set, and $\mathcal H^n(M_0)<\infty$;
            \item[(A2)'] $U \setminus M_0 = \bigcup_{i=1}^N E_{i,0}$, where $E_{i,0}$ are non-empty, open, and mutually disjoint subsets of $U$, and $N \geq 2$;
            \item[(A3)] $\partial M_0 := {\rm clos}(M_0) \setminus U$ is not empty, and for each $x \in \partial M_0$ there are at least two indexes $i_1 \neq i_2$ in $\{1,\ldots,N\}$ such that $x\in {\rm clos}\left( {\rm clos}(E_{i_j,0}) \setminus (U \cup \partial M_0) \right)$ for $j=1,2$; see Figure \ref{fig:bdry_cond}. 
        \end{itemize}
    Then, there are a Brakke flow $\mathscr{V}=\{V_t\}_{t \geq 0}$ in $U$ as well as flows $\{E_i (t)\}_{t\geq 0}$ for $i \in \{1,\ldots,N\}$ such that all conclusions of Theorem \ref{thm:canonical} hold upon setting $\Omega(x) \equiv 1$ and replacing $\R^{n+1}$ with $U$ when needed, together with the additional condition that
\begin{equation}\label{e:fixed boundary}
    {\rm clos}(\spt\|V_t\|)\setminus U = \partial M_0 \qquad \mbox{for every $t \geq 0$}\,.
\end{equation}
Moreover, $M_t \subset {\rm conv}(M_0 \cup \partial M_0)$ for every $t \geq 0$, where ${\rm conv}$ denotes convex hull.
\end{theorem}

Theorem \ref{thm:Brakke_Dirichlet} can be used to produce varifold-type solutions to Plateau's problem in $U$ (with respect to the topological spanning condition defined by \eqref{e:fixed boundary}), in a novel, purely dynamical approach, which does not rely on mass minimization nor on min-max methods.

\begin{theorem}[Asymptotic limits of Brakke flows with fixed boundary] \label{thm:asymptotics}
        Let $U$, $M_0$, and $\mathscr{V}=\{V_t\}_{t \geq 0}$ be as in Theorem \ref{thm:Brakke_Dirichlet}. There exists a sequence $\{t_j\}_{j=1}^\infty$ with $\lim_{j \to \infty} t_j = \infty$ such that $V_{t_j}$ converge as varifolds to an $n$-dimensional integral varifold $V$ in $U$ such that $V$ is stationary in $U$ and ${\rm clos}(\spt\|V\|)\setminus U = \partial M_0$. Furthermore, there are mutually disjoint open subsets $E_i \subset U$ such that $U \setminus \bigcup_{i=1}^N E_i = \spt\|V\|$, and $0 < \mathcal H^n (U \setminus \bigcup_i E_i) \leq \|V\|(U) \leq \mathcal H^n (M_0)$.
\end{theorem}

\begin{figure}[h!]
 \centering
 \includegraphics[scale=0.35]{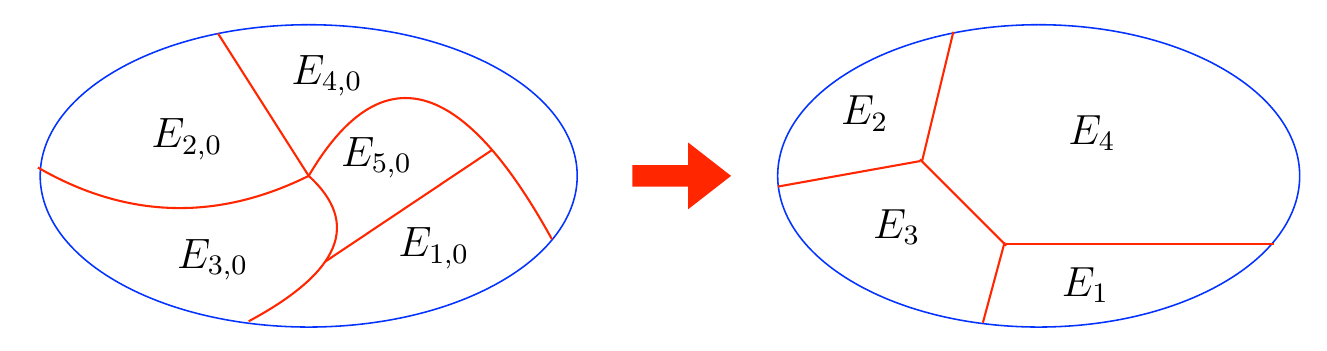}
 \caption{Asymptotic limits of Brakke flows with fixed boundary: an illustration.} \label{fig:Plateau}
\end{figure}

It is not known under which conditions one may conclude that the varifold $V$ of Theorem \ref{thm:asymptotics} is unit multiplicity. When that is the case, the flow may be potentially restarted having $V$ as an initial datum. If $V$ is a smooth minimal surface, then the regularity theory of Brakke flows together with the uniqueness theorem of smooth MCF would imply that the only flow with initial datum $V$ is the constant $V_t\equiv V$. This is not the case if $V$ is singular. In fact, we may give the following definition.
\begin{definition}[Dynamically unstable minimal surfaces] \label{def:dyn_unst}
    Let $V=\var(M,1)$ be a stationary varifold in $U$ with grain bounary structure, that is $U \setminus \spt\|V\|=\bigcup_{i=1}^N E_i$ for mutually disjoint, open subsets $E_i \subset U$. We say that $V$ is \emph{dynamically unstable} if there exists a non-trivial canonical Brakke flow in $U$ starting with $V$.
\end{definition}

The question immediately arises asking which singularity types of a minimal surface may trigger such dynamical instability. In \cite{ST_dynamical}, Tonegawa and the author proved, roughly speaking, that branching singularities always do that.

\begin{figure}
\includegraphics[scale=0.20]{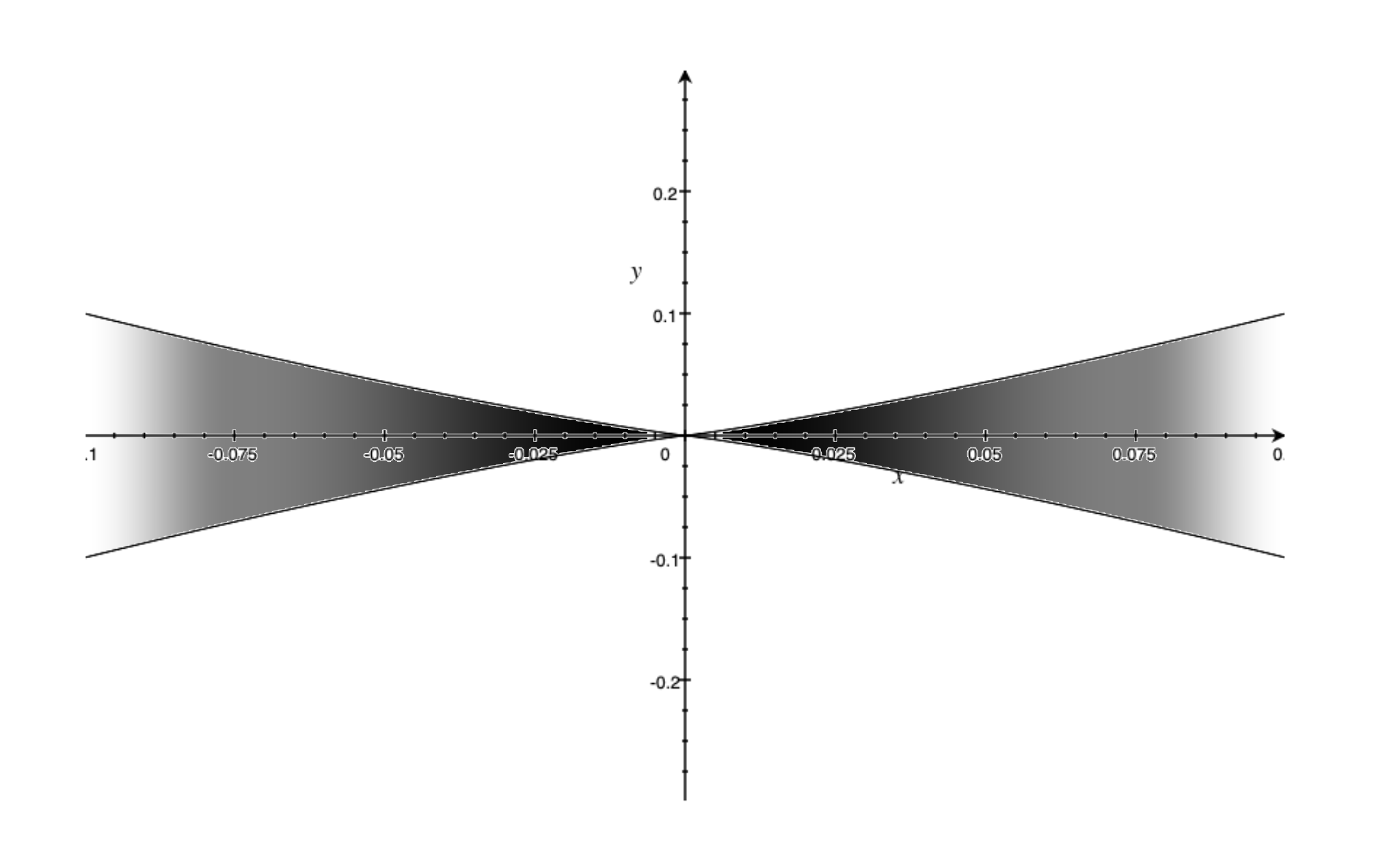} 
\caption{\small{The growth condition required in (H2): the shaded area is the region $\{\abs{x_{n+1}} \leq G(x')\}$ ($\alpha=0.51$) in the cylinder ${\rm C}_{r_0}(0,\pi_0)$ with $r_0=0.1$.}}\label{figure growth}
\end{figure}

\begin{theorem}[Dynamical instability of minimal surfaces at flat singular points] \label{thm:dynamical}
    Let $V=\var(M,1)$ be a stationary varifold in $U$ with grain boundary structure as in Definition \ref{def:dyn_unst}, and assume that there is a point $x_0 \in \spt\|V\|$, wlog $x_0=0$, such that the following holds:
    \begin{itemize}
        \item[(H1)] one of the tangent cones to $V$ at $x_0=0$ is of the form $\var(\pi_0,Q)$, for some $n$-dimensional plane $\pi_0 \in {\rm G}(n+1,n)$, wlog $\pi_0 = \R^n \times \{0\}$, and an integer $Q \ge 2$;
        \item[(H2)] there exists a radius $r_0 \in \left( 0,1 \right)$ such that, writing $x = (x',x_{n+1}) \in \R^{n+1} = \pi_0 \oplus \pi_0^\perp$ we have 
\begin{equation} \label{growth}
M \cap {\rm C}_{r_0}(0,\pi_0) \cap \{|x_{n+1}|<r_0\}\subset \{  x=(x',x_{n+1}) \in \R^{n+1} \,\colon\, \abs{x_{n+1}} \leq G(x')         \}\,,
\end{equation}
 where $G$ is the positive, radial function $G(x') = g(\abs{x'})$ defined by
\begin{equation} \label{our most affordable g}
g(s) = \frac{s}{\log^\alpha\left(1/s\right)} \qquad \mbox{for $s > 0$, with $\alpha > \frac12$}\,,
\end{equation}
see Figure \ref{figure growth}. 
    \end{itemize}
Then, there is a canonical Brakke flow $\mathscr V = \{V_t\}_{t \geq 0}$ with $\lim_{t \to 0^+} \|V_t\| = \|V_0\| = \|V\|$ and $\|V_t\|(U) < \|V\|(U)$ for every $t > 0$. In particular, $V$ is dynamically unstable.
\end{theorem}

We remark that condition ${\rm (H2)}$ implies that $\pi_0$ is the \emph{unique} tangent cone to $V$ at $x_0=0$; in fact, together ${\rm (H1)}$ and ${\rm (H2)}$ ask that the blow-ups $(\eta_r)_\sharp V$ (again with $\eta_r(x):=r^{-1}\,x$) varifold-converge to $\var(\pi_0,Q)$, as $r \to 0^+$, with rate $|\log r|^{-\alpha}$ for $\alpha > \sfrac12$. While we cannot prove that this condition is satisfied at every point of a stationary varifold admitting a flat tangent cone, there are no examples where it is violated. In fact, a faster rate of convergence of ${\rm O}(r^\alpha)$ for $\alpha \in (0,1)$ can be proved in a variety of situations, including the case of stationary varifolds that are $C^{1,\alpha}$ multiple valued minimal graphs (see \cite{SW07}), those induced by mass minimizing integral currents modulo an even integer $p \geq 4$ (see \cite{DLHMSS_decay}), and, more generally, those stationary varifolds whose regular part is stable (in the sense of non-negativity of the spectrum of the second variation operator) which satisfy ${\rm (H1)}$ at $x_0=0$ but do not have any classical singularities with density $< Q$ in a neighborhood $B_{r_0}(x_0)$ (see \cite{MWic}). In all these cases, Theorem \ref{thm:dynamical} concludes dynamical instability of the surface, and leads to speculate that those stationary varifolds $V$ for which the only canonical Brakke flow starting at $V$ is the constant flow may be free from flat singularities altogether.
 
\bibliographystyle{alpha}
\bibliography{SiGmA_refs}

\newcommand{\etalchar}[1]{$^{#1}$}
\begin{thebibliography}{DHM{\etalchar{+}}23}

\bibitem[All72]{Allard}
William~K. Allard.
\newblock On the first variation of a varifold.
\newblock {\em Ann. of Math. (2)}, 95:417--491, 1972.

\bibitem[Bra78]{Brakke_mcf}
Kenneth~A. Brakke.
\newblock {\em The motion of a surface by its mean curvature}, volume~20 of
  {\em Mathematical Notes}.
\newblock Princeton University Press, Princeton, N.J., 1978.

\bibitem[CGG91]{CGG}
Yun~Gang Chen, Yoshikazu Giga, and Shun'ichi Goto.
\newblock Uniqueness and existence of viscosity solutions of generalized mean
  curvature flow equations.
\newblock {\em J. Differential Geom.}, 33(3):749--786, (1991).

\bibitem[DGS23]{DPGS}
Guido {De Philippis}, Carlo Gasparetto, and Felix Schulze.
\newblock A short proof of {A}llard's and {B}rakke's regularity theorems.
\newblock 2023.
\newblock Preprint arXiv:2306.02490.

\bibitem[DHM{\etalchar{+}}23]{DLHMSS_decay}
Camillo {De Lellis}, Jonas Hirsch, Andrea Marchese, Luca Spolaor, and Salvatore
  Stuvard.
\newblock Excess decay for minimizing hypercurrents mod $2{Q}$.
\newblock 2023.
\newblock Preprint arXiv:2308.08704.

\bibitem[DHMS20]{DLHMS}
Camillo {De Lellis}, Jonas Hirsch, Andrea Marchese, and Salvatore Stuvard.
\newblock Regularity of area minimizing currents {${\rm mod}\,p$}.
\newblock {\em Geom. Funct. Anal.}, 30(5):1224--1336, 2020.

\bibitem[ES91]{ES}
Lawrence~C. Evans and Joel Spruck.
\newblock Motion of level sets by mean curvature. {I}.
\newblock {\em J. Differential Geom.}, 33(3):635--681, (1991).

\bibitem[ES95]{ES_IV}
Lawrence~C. Evans and Joel Spruck.
\newblock Motion of level sets by mean curvature. {IV}.
\newblock {\em J. Geom. Anal.}, 5(1):77--114, 1995.

\bibitem[HL21]{HenselLaux_ACisBV}
Sebastian Hensel and Tim Laux.
\newblock A new varifold solution concept for mean curvature flow: Convergence
  of the {A}llen-{C}ahn equation and weak-strong uniqueness.
\newblock 2021.
\newblock Preprint arXiv:2109.04233.

\bibitem[Hui90]{Huisken_mono}
Gerhard Huisken.
\newblock Asymptotic behavior for singularities of the mean curvature flow.
\newblock {\em J. Differential Geom.}, 31(1):285--299, 1990.

\bibitem[Ilm93]{Ilm_AC}
Tom Ilmanen.
\newblock Convergence of the {A}llen-{C}ahn equation to {B}rakke's motion by
  mean curvature.
\newblock {\em J. Differential Geom.}, 38(2):417--461, 1993.

\bibitem[Ilm94]{Ilm1}
Tom Ilmanen.
\newblock Elliptic regularization and partial regularity for motion by mean
  curvature.
\newblock {\em Mem. Amer. Math. Soc.}, 108(520):x+90, 1994.

\bibitem[Ilm95]{Ilmanen_mono}
Tom Ilmanen.
\newblock Singularities of mean curvature flow of surfaces.
\newblock 1995.
\newblock Preprint.

\bibitem[KT14]{Kasai-Tone}
Kota Kasai and Yoshihiro Tonegawa.
\newblock A general regularity theory for weak mean curvature flow.
\newblock {\em Calc. Var. Partial Differential Equations}, 50(1-2):1--68,
  (2014).

\bibitem[KT17]{kim-tone}
Lami Kim and Yoshihiro Tonegawa.
\newblock On the mean curvature flow of grain boundaries.
\newblock {\em Ann. Inst. Fourier (Grenoble)}, 67(1):43--142, 2017.

\bibitem[KT20]{KiTo20}
Lami Kim and Yoshihiro Tonegawa.
\newblock Existence and regularity theorems of one-dimensional {B}rakke flows.
\newblock {\em Interfaces Free Bound.}, 22(4):505--550, 2020.

\bibitem[MR08]{MugRog}
Luca Mugnai and Matthias R\"{o}ger.
\newblock The {A}llen-{C}ahn action functional in higher dimensions.
\newblock {\em Interfaces Free Bound.}, 10(1):45--78, 2008.

\bibitem[Mul56]{Mullins}
William~W. Mullins.
\newblock Two-dimensional motion of idealized grain boundaries.
\newblock {\em J. Appl. Phys.}, 27:900--904, 1956.

\bibitem[MW23]{MWic}
Paul Minter and Neshan Wickramasekera.
\newblock A structure theory for stable codimension 1 integral varifolds with
  applications to area minimising hypersurfaces mod p.
\newblock {\em J. Amer. Math. Soc.}, 2023.
\newblock Published online and ahead of print.

\bibitem[Sim83]{Simon}
Leon Simon.
\newblock {\em Lectures on geometric measure theory}, volume~3 of {\em
  Proceedings of the Centre for Mathematical Analysis, Australian National
  University}.
\newblock Australian National University, Centre for Mathematical Analysis,
  Canberra, 1983.

\bibitem[ST20]{ST_dynamical}
Salvatore Stuvard and Yoshihiro Tonegawa.
\newblock Dynamical instability of minimal surfaces at flat singular points.
\newblock 2020.
\newblock Preprint arXiv:2008.13728.

\bibitem[ST21]{ST19}
Salvatore Stuvard and Yoshihiro Tonegawa.
\newblock An existence theorem for {B}rakke flow with fixed boundary
  conditions.
\newblock {\em Calc. Var. Partial Differential Equations}, 60(1):Paper No. 43,
  53, (2021).

\bibitem[ST22]{ST_endtime}
Salvatore Stuvard and Yoshihiro Tonegawa.
\newblock End-time regularity theorem for {B}rakke flows.
\newblock 2022.
\newblock Preprint arXiv:2212.07727.

\bibitem[ST23]{ST_canonical}
Salvatore Stuvard and Yoshihiro Tonegawa.
\newblock On the existence of canonical multi-phase {B}rakke flows.
\newblock {\em Adv. Calc. Var.}, 2023.
\newblock Published online and ahead of print.

\bibitem[SW07]{SW07}
Leon Simon and Neshan Wickramasekera.
\newblock Stable branched minimal immersions with prescribed boundary.
\newblock {\em J. Differential Geom.}, 75(1):143--173, (2007).

\bibitem[Ton03]{Ton_integrality}
Yoshihiro Tonegawa.
\newblock Integrality of varifolds in the singular limit of reaction-diffusion
  equations.
\newblock {\em Hiroshima Math. J.}, 33(3):323--341, 2003.

\bibitem[Ton14]{Ton-2}
Yoshihiro Tonegawa.
\newblock A second derivative {H}\"{o}lder estimate for weak mean curvature
  flow.
\newblock {\em Adv. Calc. Var.}, 7(1):91--138, (2014).

\bibitem[Ton19]{Ton1}
Yoshihiro Tonegawa.
\newblock {\em {B}rakke's mean curvature flow: An introduction}.
\newblock SpringerBriefs in Mathematics. Springer, Singapore, 2019.

\bibitem[Whi97]{White_stratification}
Brian White.
\newblock Stratification of minimal surfaces, mean curvature flows, and
  harmonic maps.
\newblock {\em J. Reine Angew. Math.}, 488:1--35, 1997.

\bibitem[Whi05]{White_regMCF}
Brian White.
\newblock A local regularity theorem for mean curvature flow.
\newblock {\em Ann. of Math. (2)}, 161(3):1487--1519, 2005.

\end{thebibliography}
\end{document}